\input amstex
\documentstyle{amsppt}

\hsize=15.34 truecm
\vsize=22.83 truecm
\TagsOnRight
\def\d{\text{\rm d}}

\topmatter
\title{Explicit criteria for several types of ergodicity}\endtitle
\author{Mu-Fa Chen}\endauthor
\affil{(Beijing Normal University)\\
        November 7, 2000}\endaffil
\thanks {Research supported in part by 973 Project, NSFC
        and the State Education Commission of
        China.\newline
This paper is based on a talk presented at the 6th National Conference on
Probability Theory (October 14--16, 2000).}\endthanks
\address {Department of Mathematics, Beijing Normal University, Beijing 100875,
    The People's Republic of China. E-mail: mfchen\@bnu.edu.cn} \endaddress
\endkeywords
\abstract{The explicit criteria, collected in Tables 5.1 and 5.2,
for several types of ergodicity of one-dimensional diffusions or
birth-death processes have been found out recently in a
surprisingly short period. One of the criteria is for exponential
ergodicity of birth-death processes. This problem has been opened
for a long time in the study of Markov chains.  The survey article
explains in details the idea which leads to solve the problem just
mentioned. It is interesting that the problem is connected with
several branches of mathematics. Some open problems for the
further study are also proposed. }
\endabstract

\topinsert \captionwidth{10 truecm} \flushpar To appear in Chin.
J. Appl. Prob. Stat. (2001) \vskip-1.5truecm
\endinsert
\endtopmatter

\document

Let us begin with the paper by recalling
the three traditional types of ergodicity.

\subhead{1. Three traditional types of ergodicity} \endsubhead
Let $Q=(q_{ij})$ be a regular $Q$-matrix on a countable set
$E=\{i, j, k, \cdots\}$. That is, $q_{ij}\ge 0$ for all $i\ne j$,
$q_i :=-q_{ii}=\sum_{j\ne i} q_{ij}<\infty$ for all $i\in E$ and $Q$
determines uniquely a transition probability matrix $P(t)=(p_{ij}(t))$
(which is also called a $Q$-process or a Markov chain). Denote by $\pi=(\pi_i)$ a stationary distribution of $P(t)$: $\pi P(t)=\pi $ for all $t\ge 0$.
From now on, assume that the $Q$-matrix is irreducible and hence the stationary
distribution $\pi$ is unique.
Then, the three types of ergodicity are defined respectively as follows.
$$\align
&\text{\it Ordinary ergodicity}:\quad\quad \lim_{t\to  \infty}|p_{ij}(t)-\pi_j|=0 \tag 1.1\\
&\text{\it Exponential ergodicity}: \quad \lim_{t\to  \infty} e^{\hat \alpha t}|p_{ij}(t) -\pi_j| =0 \tag 1.2\\
&\text{\it Strong ergodicity}:\quad \;\;\lim_{t\to  \infty}\sup_i |p_{ij}(t)-\pi_j|=0\\
&\text{\hskip8em}\Longleftrightarrow\lim_{t\to  \infty}e^{\hat \beta t}\sup_i |p_{ij}(t)-\pi_j|=0, \tag 1.3\endalign$$
where $\hat \alpha $ and $\hat \beta $ are (the largest) positive constants and $i, j$ varies
over whole $E$. The definitions are meaningful for general Markov processes once the pointwise convergence
is replaced by the convergence in total variation norm.
The three types of ergodicity were studied in a great deal during
1953--1981. Especially, it was proved that
$\text{strong ergodicity}\Longrightarrow\text{exponential ergodicity}
\Longrightarrow \text{ordinary ergodicity}$. Refer to Anderson (1991), Chen (1992, Chapter 4)
and Meyn and Tweedie (1993) for details and
related references. The study is quite complete in the sense that we have the following criteria which
are described by the $Q$-matrix plus a test sequence $(y_i)$ only.

\proclaim{Theorem  1.1 (Criteria)} Let $H\neq\emptyset$ be an arbitrary but fixed finite subset of $E$. Then
the following conclusions hold.
\roster
\item The process $P(t)$ is ergodic iff the system of inequalities
$$\cases
\sum_j q^{}_{ij} y_j\le  -1,\quad  i\notin H\\
\sum_{i\in H}\sum_{j\ne i} q^{}_{ij} y_j < \infty
\endcases     \tag 1.4$$
has a nonnegative finite solution $(y_i)$.
\item  The process $P(t)$ is exponentially ergodic iff for some $\lambda>0$ with
$\lambda < q^{}_i$  for all $i$, the system of inequalities
$$\cases
\sum_j q^{}_{ij} y_j\le -\lambda y_i -1,\quad  i\notin H\\
\sum_{i\in H}\sum_{j\ne i} q^{}_{ij} y_j < \infty
\endcases   \tag1.5  $$
has a nonnegative finite solution $(y_i)$.
\item  The process $P(t)$ is strongly ergodic iff the system (1.4) of inequalities
has a bounded nonnegative solution $(y_i)$.
\endroster
\endproclaim

The probabilistic meaning of the criteria reads respectively as follows:
$$\max_{i\in H} \Bbb E_i \sigma_H^{}< \infty, \quad
\max_{i\in H} \Bbb E_i e^{\lambda \sigma_H^{}}< \infty
\quad\text{and}\quad \sup_{i\in E} \Bbb E_i \sigma_H^{}< \infty,$$
where $\sigma_H^{}=\inf \{t\ge \text{the first jumping time}: X_t
\in H\}$ and $\lambda$ is the same as in (1.5). The criteria are
not completely explicit since they depend on the test sequences
$(y_i)$ and in general it is often non-trivial to solve a system
of infinite inequalities. Hence, one expects to find out some
explicit criteria for some specific processes. Clearly, for this,
the first choice should be the birth-death processes. Recall that
for a birth-death process with state space $E=\Bbb Z_+=\{0, 1, 2,
\cdots\}$, its $Q$-matrix has the form: $q_{i, i+1}=b_i>0$ for all
$i\ge 0$, $q_{i, i-1}=a_i>0$ for all $i\ge 1$ and $q_{ij}=0$ for
all other $i\ne j$. Along this line, it was proved by Tweedie
(1981)(see also Anderson (1991) or Chen (1992)) that
$$ S:=\sum_{n\ge 1} \mu_n \sum_{j\le n-1} \frac{1}{\mu_j b_j}
< \infty\Longrightarrow\text{Exponential ergodicity},$$ where
$\mu_0=1$ and $\mu_n={b_0\cdots b_{n-1}}/{a_1\cdots a_n}$ for all
$n\ge 1$. Refer to Wang (1980), Yang (1986) or Hou et al (2000)
for the probabilistic meaning of $S$. However, the condition is
not necessary. A simple example is as follows. Let
$a_i=b_i=i^{\gamma }\,(i\ge 1)$ and $b_0=1$. Then the process is
exponential ergodic iff $\gamma \ge 2$ (see Chen (1996)) but $S<
\infty $ iff $\gamma
>2$. Surprisingly, the condition is correct for strong ergodicity.

\proclaim{Theorem  1.2 (H. J. Zhang, Lin and Hou (1998))}
$ S<\infty \Longleftrightarrow \,\text{Strong ergodicity}.$
\endproclaim
With a different proof, the result is extended by Y. H. Zhang (2001) to the
single-birth processes with state space $\Bbb Z_+$.
Here, the term ``single birth''means that $q_{i, i+1}>0$ for all $i\ge 0$ but
$q_{ij}\ge 0$ can be arbitrary for $j<i$.
Introducing this class of $Q$-processes is due to the
following observation:  If the first inequality in (1.4) is replaced
by equality, then we get a recursion formula for $(y_i)$ with one parameter only.
Hence, there should exist an explicit criterion for the ergodicity (resp. uniqueness, recurrence and strong ergodicity). For (1.5), there is also a recursion formula but now two parameters are involved and so it is unclear whether there exists an explicit criterion or not for the exponential ergodicity.

Note that the criteria are not enough to estimate the convergence rate $\hat \alpha$ or
$\hat \beta$ (cf. Chen (2000a)). It is the main reason why we have to come back to
study the well-developed theory of Markov chains.
For birth-death processes, the estimation of $\hat \alpha $ was studied by
Doorn in a book (1981) and in a series of papers (1985, 1987, 1991).
Especially, the precise $\hat \alpha $
was computed out for four particular models. This work has been very helpful in the later study,
as will be discussed in the next section. However, there is
no explicit criterion for $\hat \alpha >0$ ever appeared so far.
The difficulty mainly comes from the complex of two
parameters and as we will see in the next section, there is indeed an intrinsic reason for it.
By the way, we point out that there is nearly no result about the estimation of $\hat \beta $.

\proclaim{Open problem 1} How to estimate $\hat \beta$?
\endproclaim

\subhead {2. The first (non-trivial) eigenvalue (spectral
gap)}\endsubhead The birth-death processes have a nice
property---symmetriza\-bi\-lity: $\mu_i p_{ij}(t)=\mu_j p_{ji}(t)$
for all $i, j$ and $t\ge 0$.  Then, the matrix $Q$ can be regarded
as a self-adjoint operator on the real $L^2$-space $L^2(\mu)$ with
norm $\|\cdot\|$. In other words, one can use the $L^2$-theory,
which was a starting point of Doorn (1981). Furthermore, one can
study the $L^2$-exponential convergence given below. Assuming that
$\mu:=\sum_i \mu_i< \infty $ and then setting $\pi_i=\mu_i/\mu$,
we have $L^2$-space $L^2(\pi)$ with norm $\|\cdot\|$. Then, the
convergence means that
$$ \|P(t)f-\pi(f)\|\le \|f-\pi(f)\|\le e^{- \lambda_1 t} \tag 1.6$$
for all $t\ge 0$, where $\pi (f)=\int f \d \pi$ and
$ \lambda_1$ is the first non-trivial eigenvalue of $(-Q)$ (cf. Chen (1992, Chapter 9)).

The estimation of $\lambda_1$ for birth-death processes was studied by
Sullivan (1984), Liggett (1989)
and  Landim, Sethuraman and Varadhan (1996)
(see also Kipnis \& Lamdin (1999)). It was used as a comparison tool to handle the
convergence rate for some interacting particle systems, which are infinite-dimensional Markov processes.

The present author came to this topic by comparing $\lambda_1$ with $\hat \alpha $, which
is the first result in (2.1). This transfers all known results about $\hat \alpha$ to
$\lambda_1$. Then, by using the coupling methods, the author obtained
a variational formula given in the second line of (2.1).

\proclaim{Theorem  2.1}
$$\align \hat \alpha &= \lambda_1\quad\quad [\text{Chen} (1991)] \\
&=\sup_{w\in \Cal  W}\, \inf_{i\ge 0} I_i(w)^{-1}\quad
[\text{Chen} (1996)] \tag 2.1\endalign$$ where $\Cal  W =\{w:
w_i\uparrow\uparrow, \pi(w)\ge 0\}$ and $I_i(w)= \dfrac{1}{\mu_i
b_i(w_{i+1}-w_i)}\dsize \sum_{j\ge i+1} \mu_j w_j$  \endproclaim

In view of Theorem 2.1, in order to obtain a criterion for the exponential ergodicity, one needs only
to have a representative sequence $(w_i)$.
Note that the test sequence $(w_i)$ used in the formula is indeed a mimic of the eigenvector
of $\lambda_1$, and in general, the eigenvector can be very sensitive.
Thus, it becomes a question from the above formula about the existence of a representative sequence $(w_i)$ for justifying
the positiveness of $\hat \alpha$. Let us now leave Markov chains for a while and turn to diffusions.

\subhead{3. One-dimensional diffusions} \endsubhead
As a parallel of birth-death process, we now consider the elliptic operator
$L= a(x)\d^2/\d x^2+b(x)\d/\d x$ on the half line
$[0, \infty)$ with $a(x)>0$ everywhere. Again, we are interested in estimation of
the principle eigenvalues, which consist of the typical, well-known Sturm-Liouville eigenvalue
problem. Refer to  Egorov \& Kondratiev (1996) for the
present status of the study and references. Here, we mention two results, which are the most
general ones we have ever known before.

\proclaim{Theorem  3.1} Let
$b(x)\equiv 0$ (which corresponds to the birth-death process with $a_i=b_i$ for all $i\ge 1$) and set
$ \delta =\sup_{x>0} x\int_x^\infty a^{-1}$. Here we omit the integration variable when it
is integrated with respect to the Lebesgue measure. Then, we have
\roster
\item Kac \& Krein (1958):
$\delta^{-1}\ge \lambda_0 \ge (4\delta )^{-1}$,
here $\lambda_0$ is the first eigenvalue corresponding to the Dirichlet boundary $f(0)=0$.
\item Kotani \& Watanabe (1982):
$\delta^{-1}\ge \lambda_1 \ge (4\delta )^{-1}$.
\endroster
\endproclaim

In the diffusion context, we also obtained a variational formula.

\proclaim{Theorem  3.2}\,(Chen \& Wang\,(1997a), Chen(2000b)).
$$\lambda_1=\sup_{f\in \Cal  F}\, \inf_{x>0}
\bigg[\frac{e^{-C(x)}}{f'(x)}\int_x^\infty\frac{f(u)e^{C(u)}}{a(u)}\d u\bigg]^{-1}. \tag 3.1$$
where
$$\align
&C(x)=\int_0^x {b}/{a},\quad
 \pi (\d x)=\frac{1}{Z} \frac{e^{C(x)}}{a(x)}\d x \; \; \text{($Z$ is the normalizing constant) and}\\
&\Cal  F\!=\!\{f\!\in\! L^1(\pi): \pi(f)\ge 0 \text{ and }
f'|_{(0, \infty)}>0  \}.
\endalign$$
\endproclaim
In the papers quoted above, it was proved that the equality in
(3.1) replaced by ``$\ge $'' holds in general and the equality
holds whenever both $a$ and $b$ are continuous. Clearly, the
equality should also hold for certain measurable $a$ and $b$, in
virtue of a standard approximating procedure and (3.2) below.
Again, it is unclear from the result whether there exists a
representative function $f\in \Cal  F$ or not to justify the
positiveness of $\lambda_0$ or $\lambda_1$.

Before moving further, let us recall the classical variational formula
for $\lambda_1$:
$$\align \lambda_1&=\inf\{D(f):\pi(f)=0,\pi(f^2)=1\}, \tag 3.2\\
\pi(\d x)&=e^{C(x)}\d x/(a(x)Z), \quad  D(f)=\int_0^\infty a{f'}^2
\d \pi.\endalign$$ Note that there is no common point between
(3.1) and (3.2). The formula (3.1) (resp. (3.2)) is especially
meaningful in estimating the lower (resp. upper) bounds of
$\lambda_1$. Now, it is simple matter to rewrite (3.2) as (3.3)
below. Similarly, we have (3.4) for $\lambda_0$.

\proclaim {Poincar\'e inequalities}
$$\align
\lambda_1:& \quad  \text{\rm var}(f)\le \lambda_1^{-1} D(f) \tag 3.3\\
\lambda_0:& \quad \|f\|^2 \le \lambda_0^{-1} D(f), \quad  f(0)=0. \tag 3.4\endalign$$ \endproclaim

It is interesting that inequality (3.4) is a special but typical case of the weighted
Hardy inequality discussed in the next section.

\subhead{4. Weighted Hardy inequality}\endsubhead
The classical Hardy inequality goes back to Hardy (1920):
$$ \int_0^\infty \bigg(\frac{f}{x}\bigg)^p \le
\bigg(\frac{p}{p-1}\bigg)^p \int_0^\infty  {f'}^p,
\quad  f(0)= 0,f' \ge 0, $$
where the optimal constant was determined by Landau (1926). After a long
period of efforts by analysts, the inequality was finally extended to the following form.

\proclaim{Weighted Hardy inequality (Muckenhoupt (1972))}
Let $\nu $ and $\lambda $ be nonnegative Borel measures. Then
$$\int_0^\infty f^2\d \nu \le A \int_0^\infty {f'}^2 \d \lambda,
\quad  f\in C^1, f(0)=0, \tag 4.1$$
where
$B=\sup_{x>0}\nu [x,\infty]\int_x^\infty
({\d \lambda_{\text{\rm abs}}}/{\d \text{\rm Leb}})^{-1}$
satisfies $B\le A\le 4B$ and ${\d \lambda_{\text{\rm abs}}}/{\d \text{\rm Leb}}$ is
the derivative of the absolutely continuous part of $\lambda$ with respect to the Lebesgue measure.
\endproclaim

The Hardy-type inequalities play a very important role in the study of
harmonic analysis and have been treated in many publications. Refer to
the books:
Opic \& Kufner (1990),  Dynkin (1990), Mazya (1985)
and the survey article Davies (1999) for more details.
The author learnt the weighted inequality from Miclo (1999a), Bobkov
\& G\"otze (1999a, b).

It is now easy to deduce the  Poincar\'e inequality (3.4) from (4.1), simply
setting  $\nu =\pi$ and $\lambda =e^C \d x $. In other words,
Muckenhoupt's Theorem $\Longrightarrow$ Kac \& Krein's Theorem.

Thanks are given to the weighted Hardy inequality, from which we learnt that there must exist a
representative test function. Then, it is not difficult to figure out that the function is simply
$ f(x)=\sqrt{\int_0^x e^{-C}}$.

\proclaim{Theorem 4.1 (Chen (2000b, 2001a))} Let
$\delta =\sup_{x>0} \int_0^x e^{-C}\int_x^\infty {e^C}/{a}.$ Then, we have
$$\align
&(1)\;{\delta'}^{-1}  \ge \lambda_0 \ge (4\delta )^{-1}, \text{ where $\delta'$ is an explicit constant satisfying }
\delta \le \delta' \le 2 \delta. \text{ Moreover, $\lambda_0$}\\
&\quad \text{(resp. $\lambda_1$)$>0$ iff $\delta < \infty $.} \\
&(2)\;\lambda_0 =\sup_{f \in \Cal  F } \inf_{x>0} I(f)(x)^{-1}\\
& \qquad \;\; =\inf_{f \in {\Cal  F}' } \sup_{x>0} I(f)(x)^{-1}\;
\text{(Completed variational formula!).}
 \endalign$$ \endproclaim
Here, the set $\Cal  F'$ is a modification of $\Cal  F$ (The
readers are urged to find out the details from the original
papers). The lower bound in part (1) is an application of (3.1) to
the representative function mentioned above. The upper bound comes
from a direct computation. Part (1) improves (4.1) in the present
situation. Clearly, part (2) is much finer result than (4.1) and
can not be deduced from (4.1). The second line in part (2) is a
natural dual of the first one but is completely different from the
classical one (3.4).

The result can be immediately applied to the whole line or higher-dimensional situation.
For instance, for Laplacian on compact Riemannian manifolds, it was proved by Chen \& Wang\,(1997b) that
$$ \lambda_1 \ge \sup_{f \in \Cal  F}\, \inf_{r \in (0,D)} I(f)(r)^{-1} =:\xi_1$$
for some $I(f)$, which is similar to the right-hand side of (3.1).
We now have ${\delta'}^{-1} \ge \xi_1 \ge (4 \delta )^{-1}$ for some constants
$\delta $ and $\delta'$ similar to Theorem 4.1. Refer to Chen (2000b) for details.

We now return to birth-death processes. A parallel result of Theorem 4.1 has been presented in Chen (2000b).
In particular, we obtain a criterion for the positiveness of $\lambda_1$ for birth-death processes.
This is also done by Miclo (1999b) in terms of discrete Hardy inequalities. Combining this with
Theorem 2.1, we finally obtain a criterion for the exponential ergodicity of birth-death processes.
The result is included in Table 5.1 below. We mention that it is now possible to prove the last
criterion directly from (1.5), based on the study on $\lambda_1$. However, the more general case
remains open.

\proclaim {Open problem 2} Does there exist a criterion for exponential ergodicity of single-birth
processes? \endproclaim

\subhead{5. Three basic inequalities}\endsubhead
Denote by var(f)$=\|f-\pi(f)\|^2=\pi(f^2)-\pi(f)^2$. Then, the inequality (3.3)
can be rewritten as (5.1) below. On the other hand, one may study other inequalities, for instance, the logarithmic or Nash inequality listed below.
$$\align
&\text{\it Poincar\'e inequality}: \quad  \text{\rm var}(f)\le \lambda_1^{-1} D(f) \tag 5.1 \\
&\text{\it LogS inequality}: \quad  \int  f^2 \log(|f|/\|f\|)\d \pi \le \sigma^{-1} D(f) \tag 5.2\\
&\text{\it Nash inequality}: \quad  \text{\rm var}(f)^{1+2/\nu}\le
\eta^{-1} D(f)\|f \|_1^{4/\nu} \quad (\text{for some }\nu>0). \tag
5.3\endalign$$ Here, to save notation, $\sigma$ (resp. $\eta$)
denotes the largest constant so that (5.2) (resp. (5.3)) holds.

Each inequality describes a type of ergodicity. First,
$(5.1)\Longleftrightarrow(1.6)$. Next, the logarithmic Sobolev
inequality implies the decay of the semigroup $P(t)$ to $\pi$
exponentially in relative entropy with rate $\sigma$ and the Nash
inequality is equivalent to var$(P(t)(f))\le C \|f\|_1
/t^{\nu/2}$. For reversible Markov chains, a complete diagram of
these types of ergodicity was presented in Chen (1999, 2001b).

Fortunately, the criteria (also based on the weighted Hardy
inequality) for the last two inequalities as well as for the
discrete spectrum (which means that there is no continuous
spectrum and moreover, all eigenvalues have finite multiplicity)
are obtained by Mao (2000a,b,c). We can now summarize the results
in Table 5.1. The table is arranged in such order that the
property in the latter line is stranger than the former one, the
only exception is that even though the strong ergodicity is often
stronger than the logarithmic Sobolev inequality but they are not
comparable in general (Chen (2001b)).

\head Birth-death processes \endhead
$$\matrix
\format\c\;&\l\;&\l\quad&\l\;\;&\l\\
{}&i &\to i+1&\text{ at rate } & b_i=q_{i, i+1}>0\\
{}&{}&\to i-1& \text{ at rate }  &a_i=q_{i, i-1}>0.
\endmatrix$$
Define
$$\mu_0=1, \quad \mu_n=\frac{b_0\cdots b_{n-1}}{a_1\cdots a_n}, \; n\ge 1;
\quad \mu[i,k]=\sum_{i\le j\le k} \mu_j. $$
$$\vbox{\tabskip=0pt \offinterlineskip
\halign to 14.5truecm {\strut#& \vrule#\tabskip=0.5em plus 1em
&\hfil # \hfil& \vrule#& \hfil# \hfil&\vrule# \tabskip=0pt\cr
\noalign{\hrule} & & \omit\hidewidth Property \hidewidth & &
\omit\hidewidth Criterion \hidewidth &\cr\noalign{\hrule} &&
Uniqueness && $\dsize\sum_{n\ge 0}\frac{1}{\mu_n b_n} \mu[0, n] =
\infty \quad (*)$
  &\cr \noalign{\hrule}
&& Recurrence
&& $\dsize\sum_{n\ge 0} \frac{1}{\mu_n b_n}=\infty$
    &\cr \noalign{\hrule}
&& Ergodicity && $(*) \; \& \;   \mu[0, \infty)< \infty$ &\cr
\noalign{\hrule} && Exponential ergodicity && $(*) \; \& \; \dsize
\sup_{n\ge 1} \mu[n, \infty) \sum_{j\le n-1} \frac{1}{\mu_j b_j}
<\infty$
    &\cr \noalign{\hrule}
&&\! Discrete spectrum \!\!\! \! && $(*) \; \& \;
\dsize\lim_{n\to\infty}\sup_{k\ge n\!+\!1} \mu[k, \infty)
\sum_{j=n}^{k-1} \frac{1}{\mu_j b_j}=0$
    &\cr \noalign{\hrule}
&& LogS inequality && $(*) \; \& \;   \dsize\sup_{n\ge
1}\mu[n,\infty\!) \!\log[
 \mu[n,\infty\!)^{\!-\!1}]\!\!\sum_{j\le n\!-\!1}\!
    \! \frac{1}{\mu_j b_j}\!\!<\!\infty $
    &\cr \noalign{\hrule}
&& Strong ergodicity && $\!\!(*) \; \& \;   \! \dsize\sum_{n\ge
0}\!\!\frac{1}{\mu_n b_n}\! \mu[n\!\!+\!\!1,\infty\!)\!\! =\!
\dsize\!\! \sum_{n\ge 1}\!\mu_n \!\!\!\sum_{j\le n\!-\!1}\!\!
\frac{1}{\mu_j b_j}\!\!<\!\infty\!\!\!\!\!$
   &\cr \noalign{\hrule}
&& Nash inequality && $(*) \; \& \;   \dsize\sup_{n\ge 1}
\mu[n,\infty)^{\!(\nu-2)/\nu}\!\!
   \sum_{j\le n-1} \frac{1}{\mu_j b_j}\!<\!\infty (\varepsilon) $
    &\cr \noalign{\hrule}
}}$$

\centerline{Table 5.1}

\medskip

\flushpar Here, ``$(*) \; \& \;  \cdots$'' means that one requires
the uniqueness condition in the first line plus the condition
``$\cdots$''. The ``$(\varepsilon)$'' in the last line means that
there is still a small gap from being necessary.

\proclaim{Open problem 3} Find out a criterion for the  Nash
inequality in dimension one.
\endproclaim

In view of Theorem 4.1, it is natural to ask the following

\proclaim{Open question 4} Does there exist a variational formula
in the form of part $(2)$ of Theorem $4.1$ for the logarithmic
Sobolev or Nash inequality?
\endproclaim

\head Diffusion processes on $[0, \infty)$\endhead
$$L= a(x)\frac{\d^2}{\d x^2}+ b(x)\frac{\d}{\d x}.$$
Define
$$C(x)= \int_0^x b/a, \quad \quad \mu[x,y]=\int_x^y e^C/a. $$
$$\vbox{\tabskip=0pt \offinterlineskip
\halign to 14.5truecm {\strut#& \vrule#\tabskip=0.5em plus 1em
&\hfil # \hfil& \vrule#& \hfil# \hfil&\vrule# \tabskip=0pt\cr
\noalign{\hrule} & & \omit\hidewidth Property \hidewidth & &
\omit\hidewidth Criterion \hidewidth &\cr\noalign{\hrule} &&
Uniqueness && $\dsize\int_0^\infty \mu[0,x]e^{-C(x)} = \infty\quad
(*)$
  &\cr \noalign{\hrule}
&& Recurrence && $\dsize\int_0^\infty e^{-C(x)}=\infty$
    &\cr \noalign{\hrule}
&& Ergodicity && $(*) \; \& \;   \mu[0, \infty)< \infty$ &\cr
\noalign{\hrule} && Poincar\'e inequality && $(*) \; \& \;
\dsize \sup_{x> 0} \mu[x, \infty) \int_0^x e^{-C} <\infty$
    &\cr \noalign{\hrule}
&&\! Discrete spectrum \!\!\! \! && $(*) \; \& \;
\dsize\lim_{n\to\infty}\sup_{x> n} \mu[x, \infty) \int_n^x
e^{-C}=0$
    &\cr \noalign{\hrule}
&& LogS inequality && $(*) \; \& \;
\dsize\sup_{x>0}\mu[x,\infty\!) \!\log[
 \mu[x,\infty\!)^{\!-\!1}] \!\!\int_0^x e^{-C}\!<\!\infty $
    &\cr \noalign{\hrule}
&& Strong ergodicity && $\!\!(*) \; \& \;   \! \dsize\int_0^\infty
\mu[x,\infty\!)e^{-C(x)}\!<\!\infty$?
   &\cr \noalign{\hrule}
&& Nash inequality && $(*) \; \& \;   \dsize\sup_{x>0}
\mu[x,\infty)^{\!(\nu-2)/\nu}\!\!
   \int_0^x e^{-C}\!<\!\infty(\varepsilon) $
    &\cr \noalign{\hrule}
}}$$

\centerline{Table 5.2}

\medskip

\flushpar For Nash inequality, we have the same remark as before.
Note that the exponential ergodicity in Table 5.1 is replaced by
Poincar\'e inequality here. This suggests the following problem
which is conjectured to be true (partially solved in Chen (1998,
2000a). Here we mention a simple fact which was missed in the
original papers. In the reversible case, for $p_s(x,y):=\d P_s(x,
\cdot)/\d \lambda$, we have $\int(p_s(x,y)/\pi(y))^2 \pi(\d y)=
\int p_s(x,y)p_s(y,x)\lambda(\d
y)/\pi(x)=p_{2s}(x,x)/\pi(x)<\infty$. Hence $p_s(x, \cdot)/\pi\in
L^2(\pi)$).

\proclaim{Open problem 5} Prove that the exponentially ergodic convergence rate coincides with the
$L^2$-exponen\-tial convergence rate, or at least these two convergences are equivalent.\endproclaim

The question marked in the line of strong ergodicity means that it remains unproved, but conjectured to be true in terms of the similarity between one-dimensional diffusions and birth-death processes.

\proclaim{Open problem 6} Prove the criterion for strong ergodicity listed in Table 5.2.\endproclaim

\medskip

\flushpar{\bf Added in Proof}. A large part of the open problems,
except 3 and 4, have been solved by Yong-Hua Mao and Yu-Hui Zhang
in the past few months.

\enddocument